\documentclass[12pt,a4paper,twoside]{article}

\usepackage{a4wide}
\usepackage[utf8]{inputenc}
\usepackage[T1]{fontenc}
\usepackage{amsmath,amsthm,amsfonts,bbm}

\newcommand{\R}{\mathbb{R}}

\newcommand{\N}{\mathbb{N}}
\newcommand{\abs}[1]{\left\lvert{#1}\right\rvert}
\newcommand{\norm}[1]{{\left\lVert{#1}\right\rVert}}
\newcommand{\spf}[2]{{\left\langle{#1},{#2}\right\rangle}}

\newtheorem{thm}{Theorem}[section]

\begin{document}

\title{Higher-Order Space-Time Continuous Galerkin Methods for the Wave Equation}

\author{Marco Zank$^{1}$}
\date{
        $^1$Fakult\"at f\"ur Mathematik, Universit\"at Wien, \\
        Oskar-Morgenstern-Platz 1, 1090 Wien, Austria \\[1mm]
        {\tt marco.zank@univie.ac.at}
      }
      
\maketitle

\begin{abstract}
    We consider a space-time variational formulation of the second-order wave equation, where integration by parts is also applied with respect to the time variable. Conforming tensor-product finite element discretisations with piecewise polynomials of this space-time variational formulation require a CFL condition to ensure stability. To overcome this restriction in the case of piecewise multilinear, continuous ansatz and test functions, a stabilisation is well-known, which leads to an unconditionally stable space-time finite element method. In this work, we generalise this stabilisation idea from the lowest-order case to the higher-order case, i.e. to an arbitrary polynomial degree. We give numerical examples for a one-dimensional spatial domain, where the unconditional stability and optimal convergence rates in space-time norms are illustrated.
\end{abstract}

\section{Introduction}

Standard approaches for the numerical solution of hyperbolic initial-boundary value problems are usually based on semi-discretisations in space and time, where the discretisation in space and time is split accordingly. In contrast to these approaches, space-time methods discretise time-dependent partial differential equations without separating the temporal and spatial directions. In this work, the homogeneous Dirichlet problem for the second-order wave equation,
\begin{equation}\label{Zank:Einf:Gleichung}
	\left. \begin{array}{rclcl}
	\partial_{tt} u(x,t) - \Delta_xu(x,t) & = & f(x,t) & \quad & \mbox{for} \;
	(x,t) \in Q = \Omega \times (0,T), \\[1mm]
	u(x,t) & = & 0 & & \mbox{for} \; (x,t) \in \Sigma = \partial \Omega \times [0,T], 
	\\[1mm] 
	u(x,0) = \partial_t u(x,0) & = & 0 & & \mbox{for} \; x \in \Omega,
	\end{array} \right \}
\end{equation}
serves as a model problem, where $\Omega = (0,L)$ is an interval for $d=1,$ or $\Omega$ is polygonal 
for $d=2,$ or $\Omega$ is polyhedral for $d=3$, $T>0$ is a terminal time and $f$ is a given right-hand side. To derive a space-time variational formulation, we define the space-time Sobolev spaces
\begin{eqnarray*}
    H^{1,1}_{0;0,\,}(Q) &:=& L^2(0,T;H^1_0(\Omega)) \cap H^1_{0,}(0,T; L^2(\Omega)) \subset H^1(Q), \\
    H^{1,1}_{0;\,,0}(Q) &:=& L^2(0,T;H^1_0(\Omega)) \cap H^1_{,0}(0,T; L^2(\Omega))\subset H^1(Q)
\end{eqnarray*}
with the Hilbertian norms
\begin{equation*}
    \norm{v}_{H^{1,1}_{0;0,\,}(Q)} := \norm{v}_{H^{1,1}_{0;\,,0}(Q)} := \abs{v}_{H^1(Q)} := \left( \int_{0}^T \int_\Omega \big( \abs{\partial_t v(x,t)}^2 + \abs{\nabla_x v(x,t)}^2 \big) \mathrm dx \mathrm dt \right)^{1/2},
\end{equation*}
where $v \in H^1_{0,}(0,T; L^2(\Omega))$ satisfies $\|v(\cdot,0)\|_{L^2(\Omega)}=0$ and $w \in H^1_{,0}(0,T; L^2(\Omega))$ fulfils $\|w(\cdot,T)\|_{L^2(\Omega)}=0$, see \cite{ZankDissBuch2020} for more details. The bilinear form
\begin{equation*}
    a(\cdot,\cdot) \colon \, H^{1,1}_{0;0,\,}(Q) \times H^{1,1}_{0;\,,0}(Q) \to \mathbb{R},
\end{equation*}
defined by the variational identity
\begin{equation*}
	a(u,w) := -\langle \partial_t u , \partial_t w \rangle_{L^2(Q)} + \langle \nabla_x u , \nabla_x w \rangle_{L^2(Q)}
\end{equation*}
for $u \in H^{1,1}_{0;0,\,}(Q)$ and $w \in H^{1,1}_{0;\,,0}(Q)$, is continuous, i.e. the estimate
\begin{equation*}
	\forall u \in H^{1,1}_{0;0,\,}(Q) \colon \, \forall w \in H^{1,1}_{0;\,,0}(Q)\colon \quad \abs{a(u,w)} \leq \abs{u}_{H^1(Q)} \abs{w}_{H^1(Q)}
\end{equation*}
holds true due to the Cauchy-Schwarz inequality. The space-time variational formulation of \eqref{Zank:Einf:Gleichung} is to find $u \in H^{1,1}_{0;0,\,}(Q)$ such that 
\begin{equation}\label{Zank:Einf:VF}
	\forall w \in H^{1,1}_{0;\,,0}(Q)\colon \quad a(u,w) = \spf{f}{w}_{L^2(Q)},
\end{equation}
where $f \in L^2(Q)$ is a given right-hand side. Note that the initial condition 
$u(\cdot, 0) = 0$ is considered in the strong sense, whereas  the initial 
condition $\partial_t u(\cdot, 0) = 0$ is incorporated in a weak sense. The following existence and uniqueness theorem is proven in \cite[Theorem 3.2 in Chapter IV]{Ladyzhenskaya1985}, see also \cite{SteinbachZank2020Coercive, ZankDissBuch2020, Zlotnik1994}.
\begin{thm} \label{Zank:Einf:Thm:Loesbarkeit}
    For $f \in L^{2}(Q)$, a unique solution $u \in H^{1,1}_{0;0,\,}(Q)$ of the variational formulation \eqref{Zank:Einf:VF} exists and the stability estimate
    \begin{equation*}
      \abs{u}_{H^1(Q)} \leq \frac{1}{\sqrt{2}} T \norm{f}_{L^{2}(Q)}
    \end{equation*}
    holds true.
\end{thm}
Note that the solution operator
\begin{equation*}
    \mathcal L \colon \, L^2(Q) \to H^{1,1}_{0;0,\,}(Q), \quad \mathcal L f := u,
\end{equation*}
of Theorem~\ref{Zank:Einf:Thm:Loesbarkeit} is not an isomorphism, i.e. $\mathcal L$ is not surjective, see \cite{SteinbachZank2021VF, ZankDissBuch2020} for more details. In this work, for simplicity, we only consider homogeneous initial conditions, where inhomogeneous initial conditions can be treated analogously as in \cite{Ladyzhenskaya1985, Zank2021Enumath, Zlotnik1994}.

A conforming tensor-product space-time discretisation of \eqref{Zank:Einf:VF} with piecewise polynomial, continuous ansatz and test functions requires a CFL condition
\begin{equation} \label{Zank:Einf:CFL}
	h_t \leq C \, h_x
\end{equation}
with a constant $C>0$, depending on the constant of a spatial inverse inequality, where $h_t$ and $h_x$ are the mesh sizes in time and space. For a one-dimensional spatial domain $\Omega$, i.e. $d=1$, and piecewise multilinear, continuous ansatz and test functions, the CFL condition \eqref{Zank:Einf:CFL} reads as
\begin{equation*}
	h_t <  h_x
\end{equation*}
for uniform meshes with uniform mesh sizes $h_t$ and $h_x$, see \cite{SteinbachZank2020Coercive, ZankDissBuch2020}. To overcome the CFL condition \eqref{Zank:Einf:CFL}, the stabilised space-time finite element method to find $u_h \in \big( V_{h_x,0}^1(\Omega) \otimes S_{h_t}^1(0,T) \big) \cap H^{1,1}_{0;0,\,}(Q)$ such that
\begin{equation} \label{Zank:Einf:VF_disk_stab}
	-\langle \partial_t u_h , \partial_t w_h \rangle_{L^2(Q)} + \sum_{\alpha=1}^d \langle \partial_{x_\alpha} u_h , Q_{h_t}^0 \partial_{x_\alpha} w_h \rangle_{L^2(Q)} = \spf{f}{w_h}_{L^2(Q)}
\end{equation}
for all $w_h \in \big( V_{h_x,0}^1(\Omega) \otimes S_{h_t}^1(0,T) \big) \cap H^{1,1}_{0;\,,0}(Q)$ was analysed in \cite{SteinbachZankFEM2018, ZankDissBuch2020, Zlotnik1994}, where 
\begin{equation} \label{Zank:Einf:Qh0}
    Q_{h_t}^0 \colon \, L^2(Q) \to L^2(\Omega) \otimes S_{h_t}^0(0,T)
\end{equation}
is the extended $L^2$ projection on the space of the temporal piecewise constant functions and $V_{h_x,0}^1(\Omega) \otimes S_{h_t}^1(0,T)$ is the space of piecewise multilinear, continuous functions, see Section~\ref{Zank:Sec:Not} for the notations. The main results for this proposed space-time finite element method \eqref{Zank:Einf:VF_disk_stab} are the unconditional stability, i.e. no CFL condition is needed, and the space-time error estimates with
\begin{equation*}
    h:= \max\{h_x, h_t\}, \quad h_x = \max_k h_{x,k}, \quad h_t = \max_\ell h_{t,\ell},
\end{equation*}
which are summarised in the following theorem, where its proof is contained in \cite{SteinbachZankFEM2018, ZankDissBuch2020}.
\begin{thm} \label{Zank:Einf:Thm:Stab_Fehler}
  There exists a unique solution $u_h \in \big( V_{h_x,0}^1(\Omega) \otimes S_{h_t}^1(0,T) \big) \cap H^{1,1}_{0;0,\,}(Q)$ of \eqref{Zank:Einf:VF_disk_stab}, satisfying the $L^2(Q)$ stability estimate
 \begin{equation*}
    \| u_h \|_{L^2(Q)} \leq \frac{4}{\pi} T^2 \| f \|_{L^2(Q)}.
 \end{equation*}
 Further, let the solution $u$ of \eqref{Zank:Einf:Gleichung} and $\Omega$ be sufficiently regular. Then, the unique solution $u_h \in \big( V_{h_x,0}^1(\Omega) \otimes S_{h_t}^1(0,T) \big) \cap H^{1,1}_{0;0,\,}(Q)$ of \eqref{Zank:Einf:VF_disk_stab} fulfils the space-time error estimates
 \begin{align*}
    \| u - u_h \|_{L^2(Q)}  &\leq C h^2, \\
    | u - u_h |_{H^1(Q)}    &\leq C h,
 \end{align*}
 where, for the $H^1(Q)$ error estimate, a spatial inverse inequality is additionally assumed.
\end{thm}

In this work, we generalise this stabilisation idea from the linear case to the higher-order case. In greater detail, we introduce a new stabilised space-time finite element method of tensor-product type with globally continuous ansatz and test functions, which are piecewise polynomials of an arbitrary polynomial degree $p$, leading to unconditional stability and optimal convergence rates in the space-time norms $\| \cdot \|_{L^2(Q)}$, $| \cdot |_{H^1(Q)}$. In other words, the result of Theorem~\ref{Zank:Einf:Thm:Stab_Fehler} is generalised to an arbitrary polynomial degree $p$. The rest of the paper is organised as follows: In Section~\ref{Zank:Sec:Not}, notations of the used finite element spaces and $L^2$ projections are fixed. Section~\ref{Zank:Sec:FEM} introduces the new space-time finite element method. Numerical examples for a one-dimensional spatial domain and piecewise polynomials of higher-order are presented in Section~\ref{Zank:Sec:Num}. Finally, we draw some conclusions in Section~\ref{Zank:Sec:Zum}.

\section{Preliminaries} \label{Zank:Sec:Not}

In this section, notations of the used finite element spaces and $L^2$ projections are stated. For this purpose, let the bounded Lipschitz domain $\Omega \subset \R^d$ be an interval $\Omega = (0,L)$ for $d=1,$ or polygonal for $d=2,$ or polyhedral for $d=3.$ For a tensor-product ansatz, we consider admissible decompositions
\begin{equation*}
    \overline{Q} = \overline{\Omega} \times [0,T] = \bigcup_{k=1}^{N_x}\overline{\omega_k} \times \bigcup_{\ell=1}^{N_t} [t_{\ell-1},t_\ell] 
\end{equation*}
with $N:=N_x \cdot N_t$ space-time elements, where the time intervals $\tau_\ell:=(t_{\ell-1},t_\ell)$ with mesh sizes $h_{t,\ell} = t_\ell - t_{\ell-1}$ are defined via the decomposition
\begin{equation*}
    0 = t_0 < t_1 < t_2 < \dots < t_{N_t -1} < t_{N_t} = T
\end{equation*}
of the time interval $(0,T)$. The maximal and the minimal time mesh sizes are denoted by $h_t := h_{t,\max} := \max_{\ell} h_{t,\ell}$ and $h_{t,\min} := \min_{\ell} h_{t,\ell}$, respectively. For the spatial domain $\Omega$, 
we consider a shape-regular sequence $(\mathcal T_\nu)_{\nu \in {\mathbb{N}}}$ of admissible decompositions
\begin{equation*}
    \mathcal T_\nu := \{ \omega_k \subset \R^{d} \colon k=1,\dots,N_x\}
\end{equation*}
of $\Omega$ into finite elements $\omega_k \subset \R^d$ with mesh sizes $h_{x,k}$, the maximal mesh size $h_x := h_{x,\max} := \max_k h_{x,k}$ and the minimal mesh size $h_{x,\min} := \min_k h_{x,k}$. The spatial elements $\omega_k$ are 
intervals for $d=1$, triangles or quadrilaterals for $d=2$, and tetrahedra 
or hexahedra for $d=3$. 
Next, for a fixed polynomial degree $p \in \N$, we introduce the finite element space
\begin{equation*}
 Q_h^p(Q) := V_{h_x,0}^p(\Omega) \otimes S_{h_t}^p(0,T)
\end{equation*}
of piecewise polynomial, continuous functions, i.e.
\begin{equation*}
    V_{h_x,0}^p(\Omega) := V_{h_x}^p(\Omega) \cap H^1_0(\Omega) \subset H^1_0(\Omega), \quad S_{h_t}^p(0,T) \subset H^1(0,T)
\end{equation*}
with $V_{h_x}^p(\Omega) \in \left\{ S_{h_x}^p(\Omega), Q_{h_x}^p(\Omega) \right\}$. Here,
\begin{equation*}
    S_{h_t}^p(0,T) := \left\{ v_{h_t} \in C[0,T] : \forall \ell \in \{1,\dots,N_t\} \colon v_{h_t|\overline{\tau_\ell}} \in \mathbb{P}^p(\overline{\tau_\ell})   \right\}
\end{equation*}
denotes the space of piecewise polynomial, continuous functions on intervals, where $\mathbb P^p(A)$ is the space of polynomials on a subset $A \subset \R^d$ of global degree at most $p$. Analogously,
\begin{equation*}
 S_{h_x}^p(\Omega) := \left\{ v_{h_x} \in C(\overline{\Omega}) : \forall \omega \in \mathcal T_\nu \colon v_{h_x|\overline{\omega}} \in \mathbb{P}^p(\overline{\omega})   \right\}
\end{equation*}
is the space of piecewise polynomial, continuous functions on intervals ($d=1$), triangles ($d=2$), or tetrahedra ($d=3$). Moreover,
\begin{equation*}
 Q_{h_x}^p(\Omega) := \left\{ v_{h_x} \in C(\overline{\Omega}) : \forall \omega \in \mathcal T_\nu \colon v_{h_x|\overline{\omega}} \in \mathbb{Q}^p(\overline{\omega})   \right\}
\end{equation*}
is the space of piecewise polynomial, continuous functions on intervals ($d=1$), quadrilaterals ($d=2$), or hexahedra ($d=3$), where $\mathbb Q^p(A)$ is the space of polynomials on a subset $A \subset \R^d$ of degree at most $p$ in each variable. The temporal nodal basis functions of $S_{h_t}^p(0,T)$ are denoted by $\varphi_n^p$ for $n=0,\dots,p N_t$, and $\psi^p_j$, $j=1,\dots,M_x$, are the spatial nodal basis functions of $V_{h_x,0}^p(\Omega)$, i.e.
\begin{equation*}
    S_{h_t}^p(0,T) = \mbox{span} \{ \varphi_n^p \}_{n=0}^{p N_t} \quad \text{ and } \quad V_{h_x,0}^p(\Omega) = \mbox{span} \{ \psi_j^p \}_{j=1}^{M_x}.
\end{equation*}
For the stabilisation of the new space-time finite element method, we also need the spaces of piecewise polynomial, discontinuous functions
\begin{equation*}
    S^{q, \mathrm{disc}}_{h_t}(0,T) := \left\{ v_{h_t} \in L^1(0,T) :  \forall \ell \in \{1,\dots,N_t\} \colon v_{h_t|\tau_\ell} \in \mathbb{P}^q(\tau_\ell)   \right\},
\end{equation*}
where $q \in \N_0$ is a fixed polynomial degree. For a given function $v \in L^2(Q)$, the extended $L^2$ projection $Q_{h_t}^{q, \mathrm{disc}} v \in L^2(\Omega) \otimes S^{q, \mathrm{disc}}_{h_t}(0,T)$ on the space $L^2(\Omega) \otimes S^{q, \mathrm{disc}}_{h_t}(0,T)$ of piecewise polynomial, discontinuous functions with respect to the time variable is defined by
\begin{equation*}
	\spf{Q_{h_t}^{q, \mathrm{disc}} v}{v_{h_t}}_{L^2(Q)} = \spf{v}{v_{h_t}}_{L^2(Q)}
\end{equation*}
for all $v_{h_t} \in L^2(\Omega) \otimes S^{q, \mathrm{disc}}_{h_t}(0,T)$, satisfying the stability estimate
\begin{equation} \label{Zank:Not:Qh_stab}
 \| Q_{h_t}^{q, \mathrm{disc}} v \|_{L^2(Q)} \leq \| v \|_{L^2(Q)}.
\end{equation}
Note that $Q_{h_t}^0 = Q_{h_t}^{0, \mathrm{disc}}$ is the extended $L^2$ projection \eqref{Zank:Einf:Qh0} on the space of the temporal piecewise constant functions $L^2(\Omega) \otimes S^0_{h_t}(0,T) = L^2(\Omega) \otimes S^{0, \mathrm{disc}}_{h_t}(0,T)$. Analogously, for a solely time-dependent function $w \in L^2(0,T)$, we denote $Q_{h_t}^{q, \mathrm{disc}} w \in S^{q, \mathrm{disc}}_{h_t}(0,T)$ as the $L^2(0,T)$ projection on the space $S^{q, \mathrm{disc}}_{h_t}(0,T)$ of piecewise polynomial, discontinuous functions, defined by
\begin{equation*}
	\spf{Q_{h_t}^{q, \mathrm{disc}} w}{w_{h_t}}_{L^2(0,T)} = \spf{w}{w_{h_t}}_{L^2(0,T)}
\end{equation*}
for all $w_{h_t} \in S^{q, \mathrm{disc}}_{h_t}(0,T)$. We use the same notation $Q_{h_t}^{q, \mathrm{disc}}$ for solely time-dependent functions and functions, which depend on $(x,t)$, since for a function $v \in L^2(Q)$ with $v(x,t) = z(x) w(t)$, $z\in L^2(\Omega)$, $w \in L^2(0,T)$, the equality
\begin{equation*}
    Q_{h_t}^{q, \mathrm{disc}} v(x,t) = z(x) Q_{h_t}^{q, \mathrm{disc}} w(t), \quad (x,t) \in Q,
\end{equation*}
holds true.

\section{New Stabilised Space-Time Finite Element Method} \label{Zank:Sec:FEM}

In this section, we introduce a new stabilised space-time finite element method with continuous ansatz and test functions, which are piecewise polynomials of arbitrary polynomial degree $p \in \N$ with respect to the spatial variable and the temporal variable. For this purpose, we fix a polynomial degree $p \in \N$ and we introduce the perturbed bilinear form
\begin{equation*}
    a_h(\cdot,\cdot) \colon \, Q_h^p(Q) \cap H^{1,1}_{0;0,\,}(Q) \times Q_h^p(Q) \cap H^{1,1}_{0;\,,0}(Q) \to \mathbb{R}
\end{equation*}
by defining
\begin{equation*}
	a_h(u_h,w_h) := -\langle \partial_t u_h, \partial_t w_h \rangle_{L^2(Q)} + \sum_{\alpha=1}^d \langle \partial_{x_\alpha} u_h , Q_{h_t}^{p-1, \mathrm{disc}} \partial_{x_\alpha} w_h \rangle_{L^2(Q)}
\end{equation*}
for $u_h \in Q_h^p(Q) \cap H^{1,1}_{0;0,\,}(Q),$ $w_h \in Q_h^p(Q) \cap H^{1,1}_{0;\,,0}(Q)$. Note that the function $\partial_{x_\alpha} w_h$, $\alpha=1,\dots,d$, fulfils
\begin{equation*}
    \partial_{x_\alpha} w_h \in L^2(\Omega) \otimes S_{h_t}^p(0,T),
\end{equation*}
i.e. $\partial_{x_\alpha} w_h$ is still a piecewise polynomial of degree $p$ with respect to the temporal variable. The perturbed bilinear form $a_h(\cdot,\cdot)$ is continuous since the Cauchy-Schwarz inequality and the $L^2(Q)$ stability \eqref{Zank:Not:Qh_stab} of $Q_{h_t}^{p-1, \mathrm{disc}}$ yield
\begin{equation*}
    \abs{a_h(u_h,w_h)} \leq \abs{u_h}_{H^1(Q)} \abs{w_h}_{H^1(Q)}
\end{equation*}
for all $u_h \in Q_h^p(Q) \cap H^{1,1}_{0;0,\,}(Q),$ $w_h \in Q_h^p(Q) \cap H^{1,1}_{0;\,,0}(Q)$. The perturbed variational formulation, corresponding to \eqref{Zank:Einf:VF}, is to find $u_h \in Q_h^p(Q) \cap H^{1,1}_{0;0,\,}(Q)$ such that
\begin{equation} \label{Zank:FEM:VF_disk_stab_p}
    \forall w_h \in Q_h^p(Q) \cap H^{1,1}_{0;\,,0}(Q) \colon \, a_h(u_h, w_h) = \spf{f}{w_h}_{L^2(Q)}.
\end{equation}
This perturbed variational formulation \eqref{Zank:FEM:VF_disk_stab_p} coincides with the perturbed variational formulation \eqref{Zank:Einf:VF_disk_stab} for $p=1$. In other words, the new perturbed variational formulation \eqref{Zank:FEM:VF_disk_stab_p} is a generalisation of the perturbed variational formulation \eqref{Zank:Einf:VF_disk_stab} from $p=1$ to arbitrary $p \in \N.$ The numerical analysis, i.e. an analogous result as Theorem~\ref{Zank:Einf:Thm:Stab_Fehler}, of the perturbed variational formulation \eqref{Zank:FEM:VF_disk_stab_p} is far beyond the scope of this contribution, we refer to \cite{ZankWelleStabHoehereOrdnung2021}.

The discrete variational formulation \eqref{Zank:FEM:VF_disk_stab_p} is equivalent to the linear system
\begin{equation} \label{Zank:FEM:LGS}
    K_h \underline{u} = \underline{f}
\end{equation}
with the system matrix
\begin{equation*}
    K_h:= -A_{h_t} \otimes M_{h_x} + \widetilde M_{h_t} \otimes A_{h_x} \in 
{\mathbb{R}}^{M_x \cdot p N_t \times M_x \cdot p N_t},
\end{equation*}
where $M_{h_x}, \, A_{h_x} \in {\mathbb{R}}^{M_x \times M_x}$ are the mass and stiffness matrix with respect to the spatial variable, which are given by
\begin{align*}
	M_{h_x}[i,j] &= \langle \psi_j^p, \psi_i^p \rangle_{L^2(\Omega)},                     & i,j=1,\dots, M_x,  \\
	A_{h_x}[i,j] &= \langle \nabla_x \psi_j^p, \nabla_x \psi_i^p \rangle_{L^2(\Omega)},   & i,j=1,\dots, M_x,
\end{align*}
and $\widetilde M_{h_t}, \, A_{h_t} \in {\mathbb{R}}^{p N_t \times p N_t}$ 
are the perturbed mass and stiffness matrix with respect to temporal variable, which are defined by
\begin{align*}
	\widetilde M_{h_t}[n,m] &= \langle \varphi_m^p, Q_{h_t}^{p-1, \mathrm{disc}} \varphi_n^p \rangle_{L^2(0,T)}, & n=0,\dots, p N_t-1, \, m=1,\dots, p N_t, \\
	A_{h_t}[n,m] &= \langle \partial_t \varphi_m^p, \partial_t \varphi_n^p \rangle_{L^2(0,T)}, & n=0,\dots, p N_t-1, \, m=1,\dots, p N_t.
\end{align*}
Here, the nodal basis function $\varphi_0^p$ corresponds to the vertex $t_0=0$ and the nodal basis function $\varphi_{p N_t}^p$ corresponds to the vertex $t_{N_t}=T$. As the $L^2(0,T)$ projection $Q_{h_t}^{p-1, \mathrm{disc}}$ can be computed locally, i.e. on each temporal element $\tau_\ell$ for $\ell=1,\dots,N_t$, the assembling of the perturbed mass matrix $\widetilde M_{h_t}$ can be realised, as for the classical mass matrix, via local matrices.

\section{Numerical Examples} \label{Zank:Sec:Num}

In this section, numerical examples for the new space-time finite element method \eqref{Zank:FEM:VF_disk_stab_p} are given. For this purpose, we consider the hyperbolic initial-boundary value problem \eqref{Zank:Einf:Gleichung} in the one-dimensional spatial domain $\Omega := (0,1)$ with the terminal time $T=10$, i.e. the rectangular space-time
domain
\begin{equation} \label{Zank:Num:Q}
    Q:= \Omega \times(0,T) := (0,1) \times (0,10).
\end{equation}
As exact solutions, we choose
\begin{align}
    u_1(x,t) &= t^2 \sin(10 \pi x) \sin(t\,x),  \label{Zank:Num:u1}   \\
    u_2(x,t) &= t^2 (T-t)^{9/5} \sqrt{t+x^2+1} \sin (\pi  x)  \label{Zank:Num:u2}
\end{align}
for $(x,t) \in Q.$
The spatial domain $\Omega = (0,1)$ is decomposed into nonuniform elements with the vertices
\begin{equation} \label{Zank:Num:Ortsetz}
	x_0 = 0, \quad x_1 = 1/4, \quad x_2 = 1,
\end{equation}
whereas the temporal domain $(0,T) = (0,10)$ is decomposed into nonuniform elements with the vertices
\begin{equation} \label{Zank:Num:Zeinetz}
	t_0 = 0, \quad t_1 = T/8, \quad t_2 = T/4, \quad t_3 = T.
\end{equation}
We apply a uniform refinement strategy for the meshes \eqref{Zank:Num:Ortsetz}, \eqref{Zank:Num:Zeinetz}, which do not fulfil the CFL condition \eqref{Zank:Einf:CFL} at least for piecewise multilinear, continuous functions, i.e. $p=1$. Additionally, we choose $p=1$ for Table~\ref{Zank:Num:Table:u1_p1}, $p=2$ for Table~\ref{Zank:Num:Table:u1_p2}, and $p=6$ for Table~\ref{Zank:Num:Table:u1_p6}, where the number of the degrees of freedom is denoted by
\begin{equation*}
    \mathrm{dof} = M_x \cdot p \cdot N_t.
\end{equation*}
The global linear system \eqref{Zank:FEM:LGS} is solved by a direct solver, where the appearing integrals to compute the related right-hand side are calculated by using high-order quadrature rules.

In the case of piecewise multilinear, continuous functions, i.e. $p=1$, the numerical results for the smooth solution $u_1$ in \eqref{Zank:Num:u1} are given in Table~\ref{Zank:Num:Table:u1_p1}, where we observe unconditional stability, quadratic convergence in $\| \cdot \|_{L^2(Q)}$ and linear convergence in $| \cdot |_{H^1(Q)}$, as predicted by Theorem~\ref{Zank:Einf:Thm:Stab_Fehler}.

\begin{table}[!ht]
\caption{Numerical results of the Galerkin finite element discretisation \eqref{Zank:FEM:VF_disk_stab_p} for $p=1$ for the space-time cylinder \eqref{Zank:Num:Q} for the smooth function $u_1$ in \eqref{Zank:Num:u1} for a uniform refinement strategy with the starting meshes \eqref{Zank:Num:Ortsetz}, \eqref{Zank:Num:Zeinetz}.} \label{Zank:Num:Table:u1_p1}
\begin{center}
\begin{small}
\begin{tabular}{rcccccccc}
	\hline
 dof & $h_{x,\max}$ & $h_{x,\min}$ & $h_{t,\max}$ & $h_{t,\min}$ & $\norm{u_1 - u_{1,h}}_{L^2(Q)}$ & eoc & $\abs{u_1 - u_{1,h}}_{H^1(Q)}$ & eoc\\
\hline
       3 & 0.7500 & 0.2500 & 7.5000 & 1.2500 & 9.4e+01 &  -  & 2.2e+03 &  -  \\ 
      18 & 0.3750 & 0.1250 & 3.7500 & 0.6250 & 8.7e+01 & 0.1 & 2.2e+03 & 0.0 \\ 
      84 & 0.1875 & 0.0625 & 1.8750 & 0.3125 & 7.7e+01 & 0.2 & 2.0e+03 & 0.1 \\ 
     360 & 0.0938 & 0.0312 & 0.9375 & 0.1562 & 4.5e+01 & 0.8 & 1.7e+03 & 0.3 \\ 
    1488 & 0.0469 & 0.0156 & 0.4688 & 0.0781 & 1.3e+01 & 1.8 & 9.3e+02 & 0.8 \\ 
    6048 & 0.0234 & 0.0078 & 0.2344 & 0.0391 & 3.5e+00 & 1.9 & 4.9e+02 & 0.9 \\ 
   24384 & 0.0117 & 0.0039 & 0.1172 & 0.0195 & 8.8e-01 & 2.0 & 2.5e+02 & 1.0 \\ 
   97920 & 0.0059 & 0.0020 & 0.0586 & 0.0098 & 2.2e-01 & 2.0 & 1.2e+02 & 1.0 \\ 
  392448 & 0.0029 & 0.0010 & 0.0293 & 0.0049 & 5.6e-02 & 2.0 & 6.1e+01 & 1.0 \\ 
 1571328 & 0.0015 & 0.0005 & 0.0146 & 0.0024 & 1.4e-02 & 2.0 & 3.1e+01 & 1.0 \\ 
 6288384 & 0.0007 & 0.0002 & 0.0073 & 0.0012 & 3.5e-03 & 2.0 & 1.5e+01 & 1.0 \\ 
25159680 & 0.0004 & 0.0001 & 0.0037 & 0.0006 & 8.7e-04 & 2.0 & 7.7e+00 & 1.0 \\ 
\hline
\end{tabular}
\end{small}
\end{center}
\end{table}

For $p=2$ and $p=6$, the results for the smooth solution $u_1$ in \eqref{Zank:Num:u1} are stated in Table~\ref{Zank:Num:Table:u1_p2} and Table~\ref{Zank:Num:Table:u1_p6}, respectively, where we illustrate that the new space-time finite element method \eqref{Zank:FEM:VF_disk_stab_p} is unconditionally stable and the convergence rates with respect to the space-time norms $\| \cdot \|_{L^2(Q)}$, $| \cdot |_{H^1(Q)}$ are as expected. Moreover, a comparison of Table~\ref{Zank:Num:Table:u1_p1}, Table~\ref{Zank:Num:Table:u1_p2} and Table~\ref{Zank:Num:Table:u1_p6} show that a polynomial degree $p>1$ is advisable since the numbers of the degrees of freedom are much lower for $p>1$ than for $p=1$ when a fixed accuracy is desired. For example, we need $\mathrm{dof} = 25159680$ degrees of freedom for $p=1$, $\mathrm{dof} = 392448$ degrees of freedom for $p=2$ and $\mathrm{dof} = 13680$ degrees of freedom for $p=6$ to receive the error in $| \cdot |_{H^1(Q)}$ within a comparable range.

\begin{table}[!ht]
\caption{Numerical results of the Galerkin finite element discretisation \eqref{Zank:FEM:VF_disk_stab_p} for $p=2$ for the space-time cylinder \eqref{Zank:Num:Q} for the smooth function $u_1$ in \eqref{Zank:Num:u1} for a uniform refinement strategy with the starting meshes \eqref{Zank:Num:Ortsetz}, \eqref{Zank:Num:Zeinetz}.} \label{Zank:Num:Table:u1_p2}
\begin{center}
\begin{small}
\begin{tabular}{rcccccccc}
	\hline
 dof & $h_{x,\max}$ & $h_{x,\min}$ & $h_{t,\max}$ & $h_{t,\min}$ & $\norm{u_1 - u_{1,h}}_{L^2(Q)}$ & eoc & $\abs{u_1 - u_{1,h}}_{H^1(Q)}$ & eoc\\
\hline
      18 & 0.7500 & 0.2500 & 7.5000 & 1.2500 & 4.4e+03 &  -  & 1.4e+04 &  -  \\ 
      84 & 0.3750 & 0.1250 & 3.7500 & 0.6250 & 7.8e+01 & 5.8 & 2.1e+03 & 2.8 \\ 
     360 & 0.1875 & 0.0625 & 1.8750 & 0.3125 & 4.6e+01 & 0.8 & 1.7e+03 & 0.3 \\ 
    1488 & 0.0938 & 0.0312 & 0.9375 & 0.1562 & 1.2e+01 & 2.0 & 7.5e+02 & 1.2 \\ 
    6048 & 0.0469 & 0.0156 & 0.4688 & 0.0781 & 2.6e+00 & 2.2 & 2.4e+02 & 1.7 \\ 
   24384 & 0.0234 & 0.0078 & 0.2344 & 0.0391 & 2.2e-01 & 3.6 & 5.7e+01 & 2.1 \\ 
   97920 & 0.0117 & 0.0039 & 0.1172 & 0.0195 & 2.6e-02 & 3.1 & 1.4e+01 & 2.0 \\ 
  392448 & 0.0059 & 0.0020 & 0.0586 & 0.0098 & 3.2e-03 & 3.0 & 3.6e+00 & 2.0 \\ 
 1571328 & 0.0029 & 0.0010 & 0.0293 & 0.0049 & 4.0e-04 & 3.0 & 9.0e-01 & 2.0 \\ 
 6288384 & 0.0015 & 0.0005 & 0.0146 & 0.0024 & 5.1e-05 & 3.0 & 2.2e-01 & 2.0 \\ 
25159680 & 0.0007 & 0.0002 & 0.0073 & 0.0012 & 6.3e-06 & 3.0 & 5.6e-02 & 2.0 \\ 
\hline
\end{tabular}
\end{small}
\end{center}
\end{table}

\begin{table}[!ht]
\caption{Numerical results of the Galerkin finite element discretisation \eqref{Zank:FEM:VF_disk_stab_p} for $p=6$ for the space-time cylinder \eqref{Zank:Num:Q} for the smooth function $u_1$ in \eqref{Zank:Num:u1} for a uniform refinement strategy with the starting meshes \eqref{Zank:Num:Ortsetz}, \eqref{Zank:Num:Zeinetz}.} \label{Zank:Num:Table:u1_p6}
\begin{center}
\begin{small}
\begin{tabular}{rcccccccc}
	\hline
 dof & $h_{x,\max}$ & $h_{x,\min}$ & $h_{t,\max}$ & $h_{t,\min}$ & $\norm{u_1 - u_{1,h}}_{L^2(Q)}$ & eoc & $\abs{u_1 - u_{1,h}}_{H^1(Q)}$ & eoc\\
\hline
   198 & 0.7500 & 0.2500 & 7.5000 & 1.2500 & 5.2e+01 &  -  & 2.0e+03 &  -  \\ 
   828 & 0.3750 & 0.1250 & 3.7500 & 0.6250 & 3.0e+01 & 0.8 & 1.3e+03 & 0.6 \\ 
  3384 & 0.1875 & 0.0625 & 1.8750 & 0.3125 & 9.0e-01 & 5.0 & 8.6e+01 & 3.9 \\ 
 13680 & 0.0938 & 0.0312 & 0.9375 & 0.1562 & 8.9e-03 & 6.7 & 1.7e+00 & 5.6 \\ 
 55008 & 0.0469 & 0.0156 & 0.4688 & 0.0781 & 8.0e-05 & 6.8 & 3.1e-02 & 5.8 \\ 
220608 & 0.0234 & 0.0078 & 0.2344 & 0.0391 & 6.4e-07 & 7.0 & 4.9e-04 & 6.0 \\ 
883584 & 0.0117 & 0.0039 & 0.1172 & 0.0195 & 5.0e-09 & 7.0 & 7.7e-06 & 6.0 \\ 
\hline
\end{tabular}
\end{small}
\end{center}
\end{table}

For the singular solution $u_2$ in \eqref{Zank:Num:u2}, the related results are given in Table~\ref{Zank:Num:Table:u2_p1} for $p=1$, Table~\ref{Zank:Num:Table:u2_p2} for $p=2$ and  Table~\ref{Zank:Num:Table:u2_p6} for $p=6$, where we observe for $p>1$ a reduced order of convergence in  $\| \cdot \|_{L^2(Q)}$ and in $| \cdot |_{H^1(Q)}$. These convergence rates correspond to the reduced Sobolev regularity $u_2 \in H^{23/10-\varepsilon}(Q)$, $\varepsilon > 0$. 

\begin{table}[!ht]
\caption{Numerical results of the Galerkin finite element discretisation \eqref{Zank:FEM:VF_disk_stab_p} for $p=1$ for the space-time cylinder \eqref{Zank:Num:Q} for the singular function $u_2$ in \eqref{Zank:Num:u2} for a uniform refinement strategy with the starting meshes \eqref{Zank:Num:Ortsetz}, \eqref{Zank:Num:Zeinetz}.} \label{Zank:Num:Table:u2_p1}
\begin{center}
\begin{small}
\begin{tabular}{rcccccccc}
	\hline
 dof & $h_{x,\max}$ & $h_{x,\min}$ & $h_{t,\max}$ & $h_{t,\min}$ & $\norm{u_2 - u_{2,h}}_{L^2(Q)}$ & eoc & $\abs{u_2 - u_{2,h}}_{H^1(Q)}$ & eoc\\
\hline
     3 & 0.7500 & 0.2500 & 7.5000 & 1.2500 & 1.1e+03 &  -  & 4.4e+03 &  -  \\ 
    18 & 0.3750 & 0.1250 & 3.7500 & 0.6250 & 7.2e+02 & 0.6 & 2.9e+03 & 0.6 \\ 
    84 & 0.1875 & 0.0625 & 1.8750 & 0.3125 & 3.1e+02 & 1.2 & 1.4e+03 & 1.0 \\ 
   360 & 0.0938 & 0.0312 & 0.9375 & 0.1562 & 8.7e+01 & 1.8 & 5.6e+02 & 1.4 \\ 
  1488 & 0.0469 & 0.0156 & 0.4688 & 0.0781 & 2.4e+01 & 1.9 & 2.5e+02 & 1.2 \\ 
  6048 & 0.0234 & 0.0078 & 0.2344 & 0.0391 & 6.5e+00 & 1.9 & 1.1e+02 & 1.1 \\ 
 24384 & 0.0117 & 0.0039 & 0.1172 & 0.0195 & 1.6e+00 & 2.0 & 5.6e+01 & 1.0 \\ 
 97920 & 0.0059 & 0.0020 & 0.0586 & 0.0098 & 4.1e-01 & 2.0 & 2.8e+01 & 1.0 \\ 
392448 & 0.0029 & 0.0010 & 0.0293 & 0.0049 & 1.0e-01 & 2.0 & 1.4e+01 & 1.0 \\ 
1571328 & 0.0015 & 0.0005 & 0.0146 & 0.0024 & 2.6e-02 & 2.0 & 7.0e+00 & 1.0 \\ 
6288384 & 0.0007 & 0.0002 & 0.0073 & 0.0012 & 6.5e-03 & 2.0 & 3.5e+00 & 1.0 \\ 
25159680 & 0.0004 & 0.0001 & 0.0037 & 0.0006 & 1.6e-03 & 2.0 & 1.7e+00 & 1.0 \\ 
\hline
\end{tabular}
\end{small}
\end{center}
\end{table}

\begin{table}[!ht]
\caption{Numerical results of the Galerkin finite element discretisation \eqref{Zank:FEM:VF_disk_stab_p} for $p=2$ for the space-time cylinder \eqref{Zank:Num:Q} for the singular function $u_2$ in \eqref{Zank:Num:u2} for a uniform refinement strategy with the starting meshes \eqref{Zank:Num:Ortsetz}, \eqref{Zank:Num:Zeinetz}.} \label{Zank:Num:Table:u2_p2}
\begin{center}
\begin{small}
\begin{tabular}{rcccccccc}
	\hline
 dof & $h_{x,\max}$ & $h_{x,\min}$ & $h_{t,\max}$ & $h_{t,\min}$ & $\norm{u_2 - u_{2,h}}_{L^2(Q)}$ & eoc & $\abs{u_2 - u_{2,h}}_{H^1(Q)}$ & eoc\\
\hline
      18 & 0.7500 & 0.2500 & 7.5000 & 1.2500 & 5.8e+02 &  -  & 1.9e+03 &  -  \\ 
      84 & 0.3750 & 0.1250 & 3.7500 & 0.6250 & 2.0e+02 & 1.6 & 7.4e+02 & 1.4 \\ 
     360 & 0.1875 & 0.0625 & 1.8750 & 0.3125 & 3.2e+01 & 2.6 & 1.7e+02 & 2.1 \\ 
    1488 & 0.0938 & 0.0312 & 0.9375 & 0.1562 & 2.4e+00 & 3.7 & 2.7e+01 & 2.6 \\ 
    6048 & 0.0469 & 0.0156 & 0.4688 & 0.0781 & 3.9e-01 & 2.6 & 6.6e+00 & 2.0 \\ 
   24384 & 0.0234 & 0.0078 & 0.2344 & 0.0391 & 6.4e-02 & 2.6 & 2.0e+00 & 1.7 \\ 
   97920 & 0.0117 & 0.0039 & 0.1172 & 0.0195 & 1.1e-02 & 2.5 & 6.7e-01 & 1.6 \\ 
  392448 & 0.0059 & 0.0020 & 0.0586 & 0.0098 & 2.1e-03 & 2.4 & 2.4e-01 & 1.5 \\ 
 1571328 & 0.0029 & 0.0010 & 0.0293 & 0.0049 & 4.0e-04 & 2.4 & 9.3e-02 & 1.4 \\ 
 6288384 & 0.0015 & 0.0005 & 0.0146 & 0.0024 & 7.9e-05 & 2.3 & 3.7e-02 & 1.3 \\ 
25159680 & 0.0007 & 0.0002 & 0.0073 & 0.0012 & 1.6e-05 & 2.3 & 1.5e-02 & 1.3 \\ 
\hline
\end{tabular}
\end{small}
\end{center}
\end{table}

\begin{table}[!ht]
\caption{Numerical results of the Galerkin finite element discretisation \eqref{Zank:FEM:VF_disk_stab_p} for $p=6$ for the space-time cylinder \eqref{Zank:Num:Q} for the singular function $u_2$ in \eqref{Zank:Num:u2} for a uniform refinement strategy with the starting meshes \eqref{Zank:Num:Ortsetz}, \eqref{Zank:Num:Zeinetz}.} \label{Zank:Num:Table:u2_p6}
\begin{center}
\begin{small}
\begin{tabular}{rcccccccc}
	\hline
 dof & $h_{x,\max}$ & $h_{x,\min}$ & $h_{t,\max}$ & $h_{t,\min}$ & $\norm{u_2 - u_{2,h}}_{L^2(Q)}$ & eoc & $\abs{u_2 - u_{2,h}}_{H^1(Q)}$ & eoc\\
\hline
   198 & 0.7500 & 0.2500 & 7.5000 & 1.2500 & 2.7e+00 &  -  & 1.6e+01 &  -  \\ 
   828 & 0.3750 & 0.1250 & 3.7500 & 0.6250 & 6.2e-01 & 2.1 & 3.5e+00 & 2.2 \\ 
  3384 & 0.1875 & 0.0625 & 1.8750 & 0.3125 & 8.2e-02 & 2.9 & 8.8e-01 & 2.0 \\ 
 13680 & 0.0938 & 0.0312 & 0.9375 & 0.1562 & 1.5e-02 & 2.4 & 3.3e-01 & 1.4 \\ 
 55008 & 0.0469 & 0.0156 & 0.4688 & 0.0781 & 3.0e-03 & 2.3 & 1.3e-01 & 1.3 \\ 
220608 & 0.0234 & 0.0078 & 0.2344 & 0.0391 & 6.1e-04 & 2.3 & 5.3e-02 & 1.3 \\ 
883584 & 0.0117 & 0.0039 & 0.1172 & 0.0195 & 1.2e-04 & 2.3 & 2.1e-02 & 1.3 \\ 
\hline
\end{tabular}
\end{small}
\end{center}
\end{table}

\section{Conclusions} \label{Zank:Sec:Zum}

In this work, we introduced new stabilised higher-order space-time continuous Galerkin methods for the wave equation with globally continuous ansatz and test functions, which are piecewise polynomials of arbitrary polynomial degree. These methods are based on a space-time variational formulation, using also integration by parts with respect to the time variable, and its discretisation of tensor-product type with the help of a certain stabilisation. Thus, we generalised the well-known stabilisation idea from the lowest-order case to the higher-order case, i.e. to an arbitrary polynomial degree. We gave numerical examples, where the unconditional stability, i.e. no CFL condition is required, and optimal convergence rates in space-time norms were illustrated.

\bibliographystyle{acm}
\bibliography{zankliteratur}

\end{document}